\documentclass{amsart}
\usepackage{graphicx}
\usepackage{amscd}
\usepackage{amsmath}
\usepackage{amsfonts}
\usepackage{amssymb}
\newtheorem{theorem}{Theorem}
\theoremstyle{plain}

\newtheorem{definition}{Definition}
\newtheorem{example}{Example}

\newtheorem{lemma}{Lemma}

\newtheorem{proposition}{Proposition}
\newtheorem{remark}{Remark}

\numberwithin{equation}{section}

\begin{document}
\title{Local Structure of Ideal Shapes of Knots}
\author{Oguz C. Durumeric}
\address{Department of Mathematics\\
University of Iowa\\
Iowa City, Iowa 52242}
\email{odurumer@math.uiowa.edu}
\thanks{}
\date{January 21, 2002}
\dedicatory{}

\begin{abstract}
Relatively extremal knots are the relative minima of the ropelength functional
in $C^{1}$ topology. They are the relative maxima of thickness (normal
injectivity radius) functional on the set curves of fixed length, and they
include the ideal knots. We prove that a $C^{1,1}$ relatively extremal knot in
$\mathbf{R}^{n}$ either has constant maximal (generalized) curvature, or its
thickness is equal to half of the minimal double critical distance. Our main
approach is to show that the shortest curves with bounded curvature and
$C^{1\text{ }}$boundary conditions in $\mathbf{R}^{n}$ contain CLC
(circle-line-circle) curves, if they do not have constant maximal curvature. \ 
\end{abstract}\maketitle

\section{Introduction}

In this article, we study the local structure of $C^{1,1}$ relatively extremal
knots in $\mathbf{R}^{n}$ $(n\geq2),$ by using a length minimization problem
with bounded curvature and $C^{1\text{ }}$boundary conditions. The thickness
of a knotted curve is the radius of the largest tubular neighborhood around
the curve without intersections of the normal discs$.$ This is known as the
injectivity radius $i(K,\mathbf{R}^{n})$ of the normal exponential map of the
curve $K$ in the Euclidean space $\mathbf{R}^{n}$. The ideal knots are the
embeddings of $S^{1}$ into $\mathbf{R}^{n},$ maximizing $i(K,\mathbf{R}^{n})$
in a fixed isotopy (knot) class of fixed length. A relatively extremal knot is
a relative maximum of $i(K,\mathbf{R}^{n})$ in $C^{1}$ topology, if the length
is fixed.

We prove every result in $\mathbf{R}^{n}$ $(n\geq2)$ in this article, since
our methods are not dependent on dimension. However, all one dimensional knots
are trivial in $\mathbf{R}^{n}$, for $n\neq3.$ Although ideal knots are not
interesting for $n\neq3$, relatively extremal knots, the length minimization
with bounded curvature, and some of the local results on curves we obtained
may be useful for other purposes.

As noted in [Ka], ''...the average shape of knotted polymeric chains in
thermal equilibrium is closely related to the ideal representation of the
corresponding knot type''. ''Knotted DNA molecules placed in certain solutions
follow paths of random closed walks and the ideal trajectories are good
predictors of time averaged properties of knotted polymers'' as a biologist
referee pointed out to the author. Since the length of the molecule is fixed,
this problem becomes the maximization of its thickness within a fixed homotopy
class of a knot. The analytical properties of the ideal knots will be tools in
research of the physics of knotted polymers.

For simple knots, one has a good idea of the approximation of the ideal shapes
by using computers, see [Ka] and [GM]. Gonzales and Maddocks introduced the
notion of the global radius of curvature which is another characterization of
thickness in $\mathbf{R}^{3}$, and used it on discrete curves to obtain very
good approximations of the ideal shapes in [GM]. However, we do not know the
exact shape of most of the nontrivial knots (including trefoil knots) in ideal
position or the exact value of their thickness. Some estimates of the
thickness of ideal knots have been obtained by Diao [Di], Buck and Simon [BS]
and Rawdon and Simon [RS] by using results of Freedman, He and Wang [FHW].

Since a positive lower bound on thickness bounds curvature, the completion of
this class must include $C^{1,1}$ curves. The extremal cases in $\mathbf{R}%
^{3}$ are unlikely to be smooth. Very few ideal knots in $\mathbf{R}^{3}$ are
expected to be $C^{2}$, and the unknotted standard circles are possibly the
only smooth ones. This requires the study of $i(K,\mathbf{R}^{n})$ in
$C^{1,1}$ category.

In [D], the author proved the following Thickness Formula in the general
context and developed the notion of ''Geometric Focal Distance, $F_{g}(K)$''
by using metric balls, which naturally extends the notion of the focal
distance of smooth category to $C^{1}$ category. Thickness formula was first
discussed for $C^{2}-$knots in $\mathbf{R}^{3}$ in [LSDR], and for $C^{1,1}%
$-knots in $\mathbf{R}^{3}$ by Litherland in [L]. Nabutovsky, [N] has an
extensive study of $C^{1,1}$ hypersurfaces $K$ in $\mathbf{R}^{n}$ and their
injectivity radii$.$ [N] proves the upper semicontinuity of $i(K,\mathbf{R}%
^{n})$ and lower semicontinuity of $vol(K)/i(K,\mathbf{R}^{n})^{k}$ in $C^{1}
$ topology. These are also done by Litherland in [L] for $C^{1,1}-$knots in
$\mathbf{R}^{3}.$ We will use a corollary of the formula for curves in
$\mathbf{R}^{n}$ in Section 4. $i(K,M)=R_{O}(K,M),$ a rolling ball/bead
description of the injectivity radius in $\mathbf{R}^{n}$, was known by
Nabutowsky for hypersurfaces, and by Buck and Simon for $C^{2}$ curves, [BS].
The rolling ball/bead characterization is our main geometric tool. The notion
of the global radius of curvature developed by Gonzales and Maddocks for
smooth curves in $\mathbf{R}^{3}$ defined by using circles passing through 3
points of the curve in [GM] is a different characterization of $i(K,\mathbf{R}%
^{3})$ from $R_{O}$ due to positioning of the circles and metric balls.
$MDC(K)$ is the minimal double critical distance. See Section 2 for the basic definitions$.$

\bigskip

\textbf{GENERAL\ THICKNESS\ FORMULA }[D, Theorem 1]

For every complete smooth Riemannian manifold $M^{n}$ and every compact
$C^{1,1}$ submanifold $K^{k}$ $(\partial K=\emptyset)$ of $M,$
\[
i(K,M)=R_{O}(K,M)=\min\{F_{g}(K),\frac{1}{2}MDC(K)\}.
\]

\bigskip

For a $C^{1,1}$ curve $\gamma,$ $\gamma^{\prime\prime}$ exists almost
everywhere by Rademacher's Theorem, [F]. For a $C^{1,1}$ curve $\gamma(s)$
parametrized by the arclength $s$, define the (generalized) curvature
$\kappa\gamma(s)=\underset{x\neq y\rightarrow s}{\lim\sup}\frac{\measuredangle
(\gamma^{\prime}(x),\gamma^{\prime}(y))}{\mid x-y\mid}.$ $\kappa$ is defined
for all $s$. See Lemmas 1 and 2 below for a proof of $F_{g}(K^{1}%
)=F_{k}(\gamma)=\left(  \sup\text{ }\kappa\gamma\right)  ^{-1}=\left(
\sup\left\|  \gamma^{\prime\prime}\right\|  \right)  ^{-1}$ for curves
parametrized by arclength in $\mathbf{R}^{n}$, and [D, Proposition 12] for a
similar curvature description of $F_{g}(K^{k})$ for higher dimensional
$K^{k}\subset\mathbf{R}^{n}.$

Given a certain type of knot and a rope of set thickness, finding the exact
shape to tie the knot by using the shortest amount of the rope is basically
the same as finding the ideal shape of a DNA molecule of fixed length in this
knot type in $\mathbf{R}^{3}$. For any simple $C^{1,1}$ closed curve $\gamma$
in $\mathbf{R}^{n}$, define the ropelength (see [BS], [L]) or extrinsically
isoembolic length to be $\ell_{e}(\gamma)=\frac{\ell(\gamma)}{R_{o}(\gamma)}$
where $\ell(\gamma)$ is the length of $\gamma$. A curve $\gamma_{0}$ is called
an ideal (thickest) knot in a knot class $[\theta],$ if $\ell_{e}$ attains its
absolute minimum over $[\theta]\cap C^{1,1}$ at $\gamma_{0};$ and $\gamma_{0}$
is called relatively extremal, if $\ell_{e}$ attains a relative minimum at
$\gamma_{0}$ with respect to $C^{1}$ topology.

In this article, we study the pieces of relatively extremal knots away from
minimal double critical points by using minimization of length with bounded curvature.

\bigskip

\textbf{Question.} Given $p,q,v,w$ in $\mathbf{R}^{n},$ with $\left\|
v\right\|  =\left\|  w\right\|  =1$ and $\Lambda>0.$ Classify all shortest
curves in $\mathcal{C}(p,q;v,w;\Lambda)$ which is the set of all curves
$\gamma$ between the points $p$ and $q$ in $\mathbf{R}^{n}$ with
$\gamma^{\prime}(p)=v$, $\gamma^{\prime}(q)=w$ and $\kappa\gamma\leq\Lambda$.

Even though this looks like an elementary problem, a complete answer is not
known yet. This is a minimization problem with a second order differential
inequality and $C^{1}$ boundary values. There exists a shortest $C^{1,1}$
curve by Arzela-Ascoli Theorem. The following theorem classifies all cases
except the constant maximal curvature case, and also brings out the
mathematical difficulties of this problem. The results of this article are
proved by using simple geometric methods, in contrast to their analytical
nature. We include the proofs of several basic geometric facts for $C^{1,1}$ curves.

\begin{theorem}
Let $\gamma:I=[0,L]\rightarrow\mathbf{R}^{n}$ be a shortest curve in
$\mathcal{C}(p,q;v,w;\Lambda)$ parametrized by arclength.

a. If $\gamma$ does not have constant curvature $\Lambda$, i.e. $\kappa
\gamma(s_{0})<\Lambda$ for some $s_{0},$ then there exist $a_{0}$ and $b_{0}$
such that $s_{0}\in\lbrack a_{0},b_{0}]\subset\lbrack0,L]$ and $\gamma
([a_{0},b_{0}])$ is a $CLC(\Lambda)$ curve where each circular part has length
at least $\frac{\pi}{\Lambda}$ unless it contains the initial or terminal
point of $\gamma.$

b. If $R_{O}(\gamma(I),\mathbf{R}^{n})\geq\frac{1}{\Lambda}$ and $\gamma$ does
not have constant curvature $\Lambda,$ then $\gamma$ is a $CLC(\Lambda)$ curve.
\end{theorem}

A\emph{\ CLC(circle-line-circle)}$(\Lambda)$ curve is one circular arc
followed by a line segment and then by another circular arc in a $C^{1}$
fashion (like two letters J with common straight parts, one hook at each end,
and possibly non-coplanar), where the circular arcs have radius $1/\Lambda.$
If $p=q$ and $v=-w$, then the shortest curve with curvature restriction
satisfies $\kappa\equiv\Lambda$ and it is not a \emph{CLC}-curve. One can
construct curves of constant curvature $\Lambda$ with countably infinite
points where the curve is not twice differentiable. We note that the
classification of shortest curves in $\mathcal{C}(p,q;v,w;\Lambda)$ with
$\kappa\equiv\Lambda$ is not a simple matter, and it will be discussed in a
different article.

\emph{Theorem }$\emph{1}$ tells us that the parts of a relatively extremal
knot with the minimal double critical points removed are expected to be CLC
curves or overwound, i.e. $\kappa\equiv\Lambda.$ As J. Simon pointed out that
there are physical examples (no proofs) of relatively extremal unknots in
$\mathbf{R}^{3}$, which are not circles, and hence not ideal knots. One can
construct similar physical examples for composite knots.

The connectedness of the knots in \emph{Theorems 2} and \emph{3} is not
essential, and these theorems are valid on links. \emph{The General Thickness
Formula} does not assume that $K$ is connected, and \emph{Propositions 5-8} do
not use connectedness. Proof of \emph{Theorem 2} is local and based on the
nonexistence of local length decreasing and curvature nonincreasing
perturbations by repeated use of \emph{Theorem 1}. The existence of thickest
submanifolds with many components is also discussed in [D]. These proofs can
be modified by a simply changing the domain from $\mathbf{S}^{1}$ to a finite
disjoint union of circles and keeping track of which component is worked on.

\begin{theorem}
Let $\gamma:\mathbf{S}^{1}\rightarrow\mathbf{R}^{n}$ be a relatively extremal
knot, parametrized by arclength such that $\exists s_{0}\in\mathbf{S}%
^{1},\kappa\gamma(s_{0})<\sup\kappa\gamma$. Then both of the following holds.

i. $i(\gamma(\mathbf{S}^{1}),\mathbf{R}^{n})=R_{O}(\gamma(\mathbf{S}%
^{1}),\mathbf{R}^{n})=\frac{1}{2}MDC(\gamma(\mathbf{S}^{1}))$.

ii. If $s_{0}\notin I_{c}(\gamma)$, i.e. $\gamma(s_{0})$ is not a minimal
double critical point, then $\exists c\leq a<s_{0}<b\leq d$ satisfying

\qquad a) $\gamma$ is a line segment over $[a,b],$ and

\qquad b) $\gamma|[c,a]$ and $\gamma|[b,d]$ are planar circular arcs of radius
$F_{k}(\gamma)$ such that

$(a-c\geq\pi F_{k}(\gamma)$ or $c\in I_{c})$ and $(d-b\geq\pi F_{k}(\gamma)$
or $d\in I_{c}).$
\end{theorem}

As a consequence we obtain the following.

\begin{theorem}
For any relatively extremal knot $\gamma$ in $\mathbf{R}^{n}$, whose curvature
$\kappa\gamma$ is not identically $R_{O}(K)^{-1}$, the thickness of $\gamma$
is $\frac{1}{2}MDC(\gamma(\mathbf{S}^{1}))$. Equivalently, if there exists a
relatively extremal knot $\gamma$ such that $\frac{1}{2}MDC(K)>R_{O}%
(K)=F_{k}(\gamma),$ then $\gamma$ must have constant generalized curvature,
$\kappa\gamma\equiv F_{k}(\gamma)^{-1}.$
\end{theorem}

Some of our results on ideal knots overlap with [GM] which studies smooth
knots. Proposition 8 and [GM, section 4] obtain line segments away from the
maxima of the global radius of curvature $\rho_{G}$. However, maximal
$\rho_{G}$ does not distinguish between minimal double critical points and
maximal curvature points. Hence, we can obtain further conclusions, such as
\emph{Theorem 3}, and they are in a larger class ($C^{1,1}$) than smooth ideal
knots. For an ideal knot, Theorem 2 proves that (i) after a line\ segment, the
ideal curve must go through a minimal double critical point before reaching
the next line segment, and (ii) if there is a non-linear piece of the ideal
knot between a line segment and the next minimal double critical point, then
that must be a planar circular arc whose radius is the thickness of the ideal knot.

Basic definitions are given in Section 2, shortest curves with curvature
restrictions and proof of Theorem 1 are given in Section 3, and ideal knots
and proof of Theorem 2 are given in Section 4.

The author wishes to thank J. Simon and E. Rawdon for several encouraging and
helpful discussions during the completion of this work.

\section{Basic definitions for Thickness Formula}

For the generalizations of the following concepts and the Thickness Formula to
$C^{1,1}$ submanifolds of Riemannian manifolds, we refer to [D].

\begin{definition}
Let $K\subset\mathbf{R}^{n},$ be a $C^{1}$ curve in $\mathbf{R}^{n}$ and
$\gamma(s)$ be a parametrization of $K$ such that $\left\|  \gamma^{\prime
}\right\|  =1.$

i. The normal bundle to $K$ in $\mathbf{R}^{n}$ and the tangent bundle of $K$ are

$NK=\{(p,w)\in\mathbf{R}^{n}\times\mathbf{R}^{n}:s\in dom(\gamma),$
$p=\gamma(s)$ and $w\cdot\gamma^{\prime}(s)=0\}$ and

$TK=\{(p,v)\in\mathbf{R}^{n}\times\mathbf{R}^{n}:s\in dom(\gamma),$
$p=\gamma(s)$ and $v=c\gamma^{\prime}(s),c\in\mathbf{R}\},\ $respectively.

$UNK$ and $UTK$ denote all unit vectors in $NK$ and $TK,$ respectively.

$NK_{p}$ is the collection of all normal vectors of $NK$ at $p$, and the
others are defined similarly.

ii. $\exp_{p}v=p+v:$ $\mathbf{R}^{n}\times$ $\mathbf{R}^{n}\rightarrow$
$\mathbf{R}^{n}$ is the exponential map of $\mathbf{R}^{n}$ and

$\exp_{p}^{N}w=p+w:NK\rightarrow$ $\mathbf{R}^{n}$ is the normal exponential
map of $K$ into $\mathbf{R}^{n}.$
\end{definition}

\begin{definition}
i. For any metric space $X$ with a distance function $d$, $B(p,r)=\{x\in
X:d(x,p)<r\}$ and $\bar{B}(p,r)=\{x\in X:d(x,p)\leq r\}.$ For $A\subset X$ and
$x\in X$ define $d(x,A)=\inf\{d(x,a):a\in A\}$ and $B(A,r)=\{x\in
X:d(x,A)<r\}.$ The diameter $d(X)$ of $X$ is defined to be $\sup
\{d(x,y):x,y\in X\}.$ If there is ambiguity, we will use $d_{X}$ and $B(p,r;X).$

ii. For $A\subset\mathbf{R}^{n}$ and any curve $\gamma$ in $A,$ the length
$\ell(\gamma)$ is defined with respect to the metric space structure of
$\mathbf{R}^{n}.$ For any one-to-one curve $\gamma,\ell_{ab}(\gamma)$ and
$\ell_{pq}(\gamma)$ both denote the length of $\gamma$\ between $\gamma(a)=p$
and $\gamma(b)=q.$
\end{definition}

\begin{definition}
Let $K$ be a $C^{1}$ submanifold of $\mathbf{R}^{n}$. Define the thickness of
$K$ in $\mathbf{R}^{n}$ or the normal injectivity radius of $\exp^{N}$ to be
\[
i(K,M)=\sup(\left\{  0\right\}  \cup\{r>0:\exp^{N}:\{v\in NK:\left\|
v\right\|  <r\}\rightarrow M\text{ is one-to-one}\}).
\]
Equivalently, if $\gamma(s)$ parametrizes $K,$ then

$r>i(K,M)\Leftrightarrow\left(
\begin{array}
[c]{c}%
\exists\gamma(s),\gamma(t),q\in\mathbf{R}^{n},\\
\gamma(s)\neq\gamma(t),\left\|  \gamma(s)-q\right\|  <r,\left\|
\gamma(t)-q\right\|  <r,\text{ and }\\
(\gamma(s)-q)\cdot\gamma^{\prime}(s)=(\gamma(t)-q)\cdot\gamma^{\prime}(t)=0
\end{array}
\right)  .$
\end{definition}

\begin{definition}
Let $K$ be a $C^{1}$ curve in $\mathbf{R}^{n}$. For any $v\in UT\mathbf{R}%
_{p}^{n}$ and any $r>0,$ define

i. $O_{p}(v,r)=%
{\displaystyle\bigcup\limits_{w\in v^{\bot}(1)}}
B(\exp_{p}rw,r)$, where $v^{\bot}(1)=\{w\in UT\mathbf{R}_{p}^{n}:\left\langle
v,w\right\rangle =0\}$ or equivalently,

$O_{p}(v,r)=\{x\in\mathbf{R}^{n}:\exists w\in\mathbf{R}^{n},v\cdot
w=0,\left\|  w\right\|  =1,\left\|  x-p-rw\right\|  <r\}$

ii. $O_{p}^{c}(v,r)=\mathbf{R}^{n}-O_{p}(v,r)$ and

$O_{p}^{c}(v)=O_{p}^{c}(v,1)=\{x\in\mathbf{R}^{n}:\forall w\in\mathbf{R}%
^{n},v\cdot w=0,\left\|  w\right\|  =1,\left\|  x-p-w\right\|  \geq1\}$

iii. $O_{p}(r;K)=O_{p}(v,r)$ where $v\in UTK_{p}$

iv. $O(r;K)=%
{\displaystyle\bigcup\limits_{p\in K}}
O_{p}(r;K)$
\end{definition}

\begin{definition}
Let $K$ be a $C^{1}$ curve in $\mathbf{R}^{n}.$ Define

i. The ball radius of $K$ in $\mathbf{R}^{n}$ to be $R_{O}(K,\mathbf{R}%
^{n})=\inf\{r>0:O(r;K)\cap K\neq\emptyset\}$

ii. The pointwise geometric focal distance $F_{g}(p)=\inf\{r>0:p\in
\overline{O_{p}(r;K)\cap K}\}$ for any $p\in K,$ and the geometric focal
distance $F_{g}(K)=\inf_{p\in K}F_{g}(p).$
\end{definition}

\begin{definition}
Let $K$ be a $C^{1}$ submanifold of $M.$ A pair of points $p$ and $q$ in $K$
are called a double critical pair for $K,$ if there is a line segment
$\gamma_{pq}$ of positive length between $p$ to $q$, normal to $K$ at both $p
$ and $q.$ Define the minimal double critical distance

$MDC(K)=\inf\{\left\|  p-q\right\|  :\{p,q\}$ is a double critical pair for $K\}.$
\end{definition}

\section{Shortest Curves in $\mathbf{R}^{n}$ with Curvature Restrictions}

In this section, $\gamma:I\rightarrow\mathbf{R}^{n}$ denotes a simple $C^{1}$
curve with $I=[0,L]$, $\left\|  \gamma^{\prime}\right\|  \neq0$ and
$K=image(\gamma)$.

\begin{definition}
Let $\{e_{i}:i=1,2,..n\}$ be the standard basis in $\mathbf{R}^{n}$. Let
$E_{i},$ $E_{i}^{+},$ $E_{i}^{-}$ denote the $e_{i}-axis,$ its positive and
negative parts.
\end{definition}

\begin{definition}
For $\gamma:I\rightarrow\mathbf{R}^{n}$, define:

Dilations: $dil^{d}\gamma^{\prime}(s,t)=\frac{\left\|  \gamma^{\prime
}(s)-\gamma^{\prime}(t)\right\|  }{\ell_{st}(\gamma)}$ and $dil^{\alpha}%
\gamma^{\prime}(s,t)=\frac{\measuredangle(\gamma^{\prime}(s),\gamma^{\prime
}(t))}{\ell_{st}(\gamma)}$ for $s\neq t$

Curvature: $\kappa\gamma(s)=\underset{t,u\rightarrow s}{\lim\sup\text{ }%
}dil^{\alpha}\gamma^{\prime}(t,u)$

Lower curvature: $\kappa^{-}\gamma(s)=\underset{t\rightarrow s}{\lim\sup\text{
}}dil^{\alpha}\gamma^{\prime}(s,t)$

Analytic focal distance: $F_{k}(\gamma)=(\underset{}{\sup_{I}}\kappa
\gamma(s))^{-1}$

\begin{remark}
i. Since $\underset{v\rightarrow w}{\lim}\frac{\measuredangle(v,w)}{\left\|
v-w\right\|  }=1$ for $\left\|  v\right\|  =\left\|  w\right\|  =1,$ one
obtains the same $\kappa\gamma$, if one uses $dil^{d}$ instead of $dil^{a}$,
provided that $\gamma$ is parametrized with respect to the arclength. Same is
true for $\kappa^{-}\gamma.$

ii. $\kappa\gamma(s)\geq\kappa^{-}\gamma(s),\forall s.$

iii. If $\gamma\in C^{1,1}$ and $\left\|  \gamma^{\prime}\right\|  \equiv1, $
then $\left\|  \gamma^{\prime\prime}(s)\right\|  =\kappa^{-}\gamma(s),\forall
s $ a.e.

iv. $\underset{s_{n}\rightarrow s}{\lim\sup}$ $\kappa\gamma(s_{n})\leq
\kappa\gamma(s)$.
\end{remark}
\end{definition}

\begin{lemma}
All of the following are equivalent for a $C^{1}$ curve $\gamma:I\rightarrow
\mathbf{R}^{n}$ with $\left\|  \gamma^{\prime}\right\|  \equiv1,$ for the same
$\Lambda.$

i. $\kappa\gamma(s)\leq\Lambda,\forall s\in I.$

ii. $dil^{d}\gamma^{\prime}(s,t)\leq\Lambda,\forall s,t\in I.$

iii. $dil^{\alpha}\gamma^{\prime}(s,t)\leq\Lambda,\forall s,t\in I.$

iv. $\left\|  \gamma^{\prime\prime}(s)\right\|  \leq\Lambda,$ $\forall s\in I$
a.e., and $\gamma^{\prime}$is absolutely continuous.
\end{lemma}

\begin{proof}
$(i\Longrightarrow iii):$ $\forall s<t<u,$ $dil^{\alpha}\gamma^{\prime
}(s,u)\leq\max(dil^{\alpha}\gamma^{\prime}(s,t),dil^{\alpha}\gamma^{\prime
}(t,u)).$ Hence, if $dil^{\alpha}\gamma^{\prime}(s,u)\geq A$ for some $s\neq
u$ and $A,$ then there exists $s_{0}\in\lbrack s,u]$ with $\kappa\gamma
(s_{0})\geq A.$

$(iv\Longrightarrow ii):\left\|  \gamma^{\prime}(t)-\gamma^{\prime
}(s)\right\|  =\left\|  \int\limits_{s}^{t}\gamma^{\prime\prime}(u)du\right\|
\leq\int\limits_{s}^{t}\left\|  \gamma^{\prime\prime}(u)\right\|
du\leq\Lambda\left\|  t-s\right\|  ,$ by absolute continuity.

$(iii\Longrightarrow i)$ and $(ii\Longrightarrow iv)$ are obvious, and
$(ii\Longleftrightarrow iii)$ is by Remark 1.i.
\end{proof}

\begin{definition}
A curve $\gamma:I\rightarrow\mathbf{R}^{n}\ $is called to have curvature at
most $\Lambda,$ if $\kappa\gamma\leq\Lambda$ on $I.$ By the previous lemma,
$\gamma$ must be of class $C^{1,1}.$
\end{definition}

\begin{definition}
A $C^{1,1}$ curve $\gamma:I=[a,b]\rightarrow\mathbf{R}^{n}\ $is called \ a
$CLC(\Lambda)$ curve if there are $a\leq c\leq d\leq b$ such that (i)
$\gamma([c,d])$ is a line segment of possibly zero length, and (ii) each of
$\gamma([a,c])$ and $\gamma([d,b])$ is a planar circular arc of radius
$\frac{1}{\Lambda}$ and of length in $[0,\frac{2\pi}{\Lambda}),$ with the
possibility that $\gamma$ is not planar.
\end{definition}

\begin{proposition}
Let $p\in\mathbf{R}^{n}$ and $v\in\mathbf{R}^{n}$ with $\left\|  v\right\|
=1. $

a) $\forall q\in O_{p}^{c}(v),\exists w\in\mathbf{R}^{n},\exists q^{\prime}%
\in\partial O_{p}^{c}(v)$ and a $C^{1}$ curve $\gamma_{pq}\subset span\{v,w\}$
such that

i. $v\cdot w=0,$ $\left\|  w\right\|  =1,$ and

ii. $\gamma_{pq}(t)=\left\{
\begin{array}
[c]{cc}%
v\sin t+w(1-\cos t) & \text{if }0\leq t\leq t_{0}\\
q^{\prime}+(t-t_{0})(q-q^{\prime}) & \text{if }t_{0}\leq t\leq t_{1}%
\end{array}
\right.  ,$ where $q^{\prime}=\gamma_{pq}(t_{0})$ and $q=\gamma_{pq}(t_{1}),$ and

iii. $\gamma_{pq}$ is a shortest curve among all the continuous curves
$\varphi$ from $p$ to $q$ in $O_{p}^{c}(v)$ with $(\varphi(t)-p)\cdot v>0$ for
small $t>0$ and $\varphi(0)=p$.

b) If $q-p\neq\lambda v,\forall\lambda,$ then $q^{\prime}$ and $w$ are unique
and $\gamma_{pq}$ is unique up to parametrization.

\begin{proof}
Consider a non-empty set of all rectifiable curves of length $\leq L$
satisfying (a.iii). Parametrize each curve $\varphi$ by arclength and extend
the domain to $[0,L]$ by keeping $\varphi$\ constant after reaching $q$ so
that $\left\|  \varphi(s)-\varphi(t)\right\|  \leq\left|  s-t\right|  ,\forall
s,t.$ This forms a non-empty, bounded and equicontinuous family, and length
functional is lower semi-continuous under uniform convergence. By
Arzela-Ascoli Theorem, a shortest curve $\gamma_{pq}$ from $p$ to $q$ in
$O_{p}^{c}(v)$ satisfying (a.ii) exists. Also, the proof below shows how to
deform any curve $\varphi$ as in (a.iii) to a shorter curve, where the aim is
to reach $\gamma_{pq}.$

It suffices to give the rest of the proof for $p=0,$ $v=e_{1}$ and $q\neq0.$
Set $\gamma_{pq}=\gamma.$

\textbf{Case 1.} If $q=\lambda e_{1}$, $\lambda>0,$ then $\gamma$ is the line
segment from $0$ to $q$ where $q^{\prime}=0$ and $t_{0}=0$. Conversely, if
$\gamma$ intersects $E_{1}^{+}$ at $q^{\prime\prime}\neq0$, then $q=\lambda
e_{1}$ for some $\lambda>0.$ For, $\gamma$ must be along $E_{1}^{+} $ between
$0$ and $q^{\prime\prime}$, and then extends uniquely as a geodesic of
$\mathbf{R}^{n}$ beyond $q^{\prime\prime}.$

For any $u\in\mathbf{R}^{n}$, define $u^{N}=u-(u\cdot e_{1})e_{1}.$

\textbf{Case 2. }$\gamma\cap E_{1}=\{0\}.$ Thus, $q^{N}\neq0$ and define
$w=q^{N}/\left\|  q^{N}\right\|  .$ It suffices to give the proof for
$w=e_{2}.$ Define $f:\mathbf{R}^{n}-E_{1}\rightarrow A=\{xe_{1}+ye_{2}%
:x,y\in\mathbf{R} $ and $y>0\}$ by $f(u)=(u\cdot e_{1})e_{1}+\left\|
u^{N}\right\|  e_{2}.$ $f$ is a length decreasing map:
\begin{align*}
\left\|  f(u)-f(z)\right\|  ^{2}  &  =\left\|  (u\cdot e_{1})e_{1}-(z\cdot
e_{1})e_{1}\right\|  ^{2}+\left(  \left\|  u^{N}\right\|  -\left\|
z^{N}\right\|  \right)  ^{2}\\
&  \leq\left\|  (u-z)\cdot e_{1}\right\|  ^{2}+\left\|  u^{N}-z^{N}\right\|
^{2}=\left\|  u-z\right\|  ^{2}%
\end{align*}
and equality holds if and only if $u^{N}=cz^{N}$ for some $c>0,$ i.e. $u\in
span(e_{1},z).$ Reparametrize $\gamma$ with respect to arclength.\ By
following Federer [F], pp. 109, 163-168, we obtain that $\gamma$ is lipschitz,
absolutely continuous, $\gamma^{\prime}$ exists a.e. and
\[
\ell(\gamma)=\int\limits_{0}^{\ell(\gamma)}\left\|  \gamma^{\prime
}(s)\right\|  ds\geq\int\limits_{0}^{\ell(\gamma)}\left\|  f_{\ast}%
\gamma^{\prime}(s)\right\|  ds\geq\ell(f(\gamma)).
\]
Since $\gamma$ is a shortest curve from $0$ to $q$, $\left\|  \gamma^{\prime
}(s)\right\|  =\left\|  f_{\ast}\gamma^{\prime}(s)\right\|  $ and
$\gamma^{\prime}(s)\in span\{e_{1},\gamma(s)\}$ for almost all $s\in
\lbrack0,\ell(\gamma)].$
\[
(\gamma^{N})^{\prime}(s)=\gamma^{\prime}(s)^{N}=\lambda(s)\gamma^{N}(s),\text{
for }s\in\lbrack0,\ell(\gamma)],a.e.
\]%
\[
\frac{d}{ds}(\gamma^{N}(s)(\gamma^{N}(s)\cdot\gamma^{N}(s))^{-\frac{1}{2}%
})=0,\text{ for }s\in\lbrack0,\ell(\gamma)],a.e.
\]
By absolute continuity and $\gamma(\ell(\gamma))=q\in span\{e_{1},e_{2}\}, $
one obtains that $\gamma(I)\subset span\{e_{1},e_{2}\}$. This reduces the
proof to the $\mathbf{R}^{2}$ case.

\textbf{Subcase 2.1. }$\left\|  q-e_{2}\right\|  =1$, that is $q\in\partial
O_{p}^{c}(v).$ Define $g:\{u\in\mathbf{R}^{2}:\left\|  u-e_{2}\right\|
\geq1\}\rightarrow\{u\in\mathbf{R}^{2}:\left\|  u-e_{2}\right\|  =1\}$ by
$g(u)=e_{2}+\frac{u-e_{2}}{\left\|  u-e_{2}\right\|  }.$ Then, $g$ is a
distance decreasing map, $\left\|  g(u)-g(z)\right\|  \leq\left\|
u-z\right\|  $ and equality holds if and only if $\left\|  u-e_{2}\right\|
=\left\|  z-e_{2}\right\|  =1.$ Hence, $\ell(\gamma)\geq\ell(g(\gamma))$, and
consequently the shortest curve $\gamma$\ must lie on the circle $\left\|
u-e_{2}\right\|  =1$ between $p$ and $q$, by a proof similar to above with $f$.

\textbf{Subcase 2.2.} $\left\|  q-e_{2}\right\|  >1$, that is $q\in
intO_{p}^{c}(v).$ Any component of $\gamma\cap intO_{p}^{c}(v)$ is a line
segment. Let $\eta$ be the component containing $q$. By case assumption and
Case 1, $\bar{\eta}\cap E_{1}^{+}=\emptyset.$ There exists unique $q\prime$ in
$\bar{\eta}\cap\partial O_{p}^{c}(v)$ with $\left\|  q^{\prime}-e_{2}\right\|
=1.$ By Case 2.1, $\gamma$ is a union of a segment and a circular arc. If
$\gamma$ were not $C^{1}$ at $q\prime=\gamma(t_{0})$, then for sufficiently
small $\varepsilon>0,$ the line segments between $\gamma(t_{0}-\varepsilon)$
and $\gamma(t_{0}+\varepsilon)$ lie in $O_{p}^{c}(v)$ and have length
$<2\varepsilon$, by the first variation. Hence, $\gamma$ is $C^{1},$ satisfies
a.i-iii, in $\mathbf{R}^{2}=span\{e_{1},e_{2}\}$ and consequently in
$\mathbf{R}^{n}$. In Case 2, $q^{\prime}$ and $w$ are unique and $\gamma_{pq}$
is unique up to parametrization.

\textbf{Case 3.} $\gamma\cap E_{1}^{-}\neq\emptyset.$ \textbf{Subcase 3.1.
}$q\in E_{1}^{-}.$ $q\in intO_{p}^{c}(v)$ and let $\eta$ be the line segment
part of $\gamma$ ending at $q.$ Obviously, $\eta\nsubseteq E_{1}^{-}.$ Choose
$q^{\prime\prime}\in(\eta-\{q\}-\partial O_{p}^{c}(v)).$ $\gamma$ restricts
the shortest curve from $p$ to $q^{\prime\prime}$, by Case 2. Hence, $\gamma$
\ follows a circular arc to $q^{\prime}$ then a segment to $q^{\prime\prime}%
,$which must be $\eta.$ This proves (a.i-iii). By rotating $\gamma$\ around
$E_{1}$, one obtains infinitely many shortest curves $\gamma_{\alpha}$
satisfying (a.i-iii).

\textbf{Subcase 3.2.} Suppose there exists $q^{\prime\prime\prime}\in
\gamma\cap E_{1}^{-}$ and $q^{\prime\prime\prime}\neq q.$ Then by following
any $\gamma_{\alpha}\neq\gamma$ from $p$ to $q^{\prime\prime\prime},$ and
$\gamma$ from $q^{\prime\prime\prime}$ to $q$, creates a shortest curve with a
corner within $intO_{p}^{c}(v),$ an open subset of $\mathbf{R}^{n}.$ Hence,
Subcase 3.2 does not occur.
\end{proof}
\end{proposition}

\begin{proposition}
Let $\gamma:I=[0,L]\rightarrow\mathbf{R}^{n}$ be with $\kappa\gamma\leq1 $ and
$\left\|  \gamma^{\prime}\right\|  \equiv1.$ Then,

a. $\gamma(s)\in O_{\gamma(a)}^{c}(\gamma^{\prime}(a))$, $\forall a,s\in I$
with $\left|  s-a\right|  \leq\pi.$ Also,

$(\gamma(s_{0})\in\partial O_{\gamma(a)}^{c}(\gamma^{\prime}(a))$ for some
$a,s_{0}\in I$ with $0<\left|  s_{0}-a\right|  \leq\pi)$ if and only if

$\gamma$ is a circular arc of radius $1$ in $\partial O_{\gamma(a)}^{c}%
(\gamma^{\prime}(a))$ between $\gamma(a)$ and $\gamma(s_{0}).$

b. If $\left\|  \gamma(0)\right\|  =\left\|  \gamma(L)\right\|  =1$ and
$\left\|  \gamma(a)\right\|  >1$ for some $a\in\lbrack0,L]$, then $L>\pi.$

c. If $\gamma^{\prime\prime}(a)$ exists and $\left\|  \gamma^{\prime\prime
}(a)\right\|  =1,$ for some $a\in\lbrack0,L)$ then

$\qquad\forall R>1,\exists\varepsilon>0$ such that $\gamma((a,a+\varepsilon
))\subset B(\gamma(a)+R\gamma^{\prime\prime}(a),R).$

\begin{proof}
The proof follows the following order: (a) for $0\leq\left|  s-a\right|
\leq\frac{\pi}{2},$ (b) is next, and then (a) for $\left|  s-a\right|  \leq
\pi.$ (c) is independent.

\textbf{(a:}$\frac{\pi}{2}$\textbf{)} By using an isometry of $\mathbf{R}^{n}
$, reparametrization and symmetry, it suffices to prove this for $a=0,$
$\gamma(0)=0$, $\gamma^{\prime}(0)=e_{1}$ and for $0\leq s\leq\frac{\pi}{2}.$

$\gamma^{\prime}(s)=\alpha(s)e_{1}+\beta(s)v(s)$ where $\left\|  v(s)\right\|
=1$ and $v(s)\cdot e_{1}=0,$ for $s\in\lbrack0,\frac{\pi}{2}].$ Then, by Lemma
1, $\measuredangle(\gamma^{\prime}(0),\gamma^{\prime}(s))\leq s,$
$\alpha(s)\geq\cos s,$ and $\beta(s)\leq\sin s,$ since $\alpha^{2}+\beta
^{2}=1.$ For any $u\in\mathbf{R}^{n}$, define $u^{N}=u-(u\cdot e_{1})e_{1}.$%
\begin{align*}
\gamma(s)\cdot e_{1}  &  =\int\nolimits_{0}^{s}\gamma^{\prime}(t)\cdot
e_{1}dt\geq\sin s\\
\left\|  \gamma(s)^{N}\right\|   &  \leq\int\nolimits_{0}^{s}\left|
\beta(t)\right|  dt\leq1-\cos s
\end{align*}
For any unit vector $u$ normal to $e_{1},$

$\left\|  u-\gamma(s)\right\|  ^{2}=\left(  \gamma(s)\cdot e_{1}\right)
^{2}+\left\|  u-\gamma(s)^{N}\right\|  ^{2}\geq\sin^{2}s+\left\|  u-\left\|
\gamma(s)^{N}\right\|  u\right\|  ^{2}\geq1.$ Hence, $\gamma(s)\in O_{0}%
^{c}(e_{1})$ for $0\leq s\leq\frac{\pi}{2}.$ Suppose that $\gamma(s_{0}%
)\in\partial O_{0}^{c}(e_{1})$ for some $s_{0}\in(0,\frac{\pi}{2}].$ Then, all
of the above inequalities become equalities for a fixed $u$ and $\gamma
(s)^{N}$ is parallel to $u,$ to conclude $\gamma(s)=\left(  \sin s\right)
e_{1}+(1-\cos s)u,$ for $s\in(0,s_{0}]$.

\textbf{(b) }Choose $m\in\lbrack0,L],$ such that $\left\|  \gamma(m)\right\|
\geq\left\|  \gamma(s)\right\|  ,\forall s\in\lbrack0,L].$ Since $\gamma(m)$
is a furthest point from $\mathbf{0},$ $\gamma^{\prime}(m)\cdot\gamma
(m)=0\ $and $\mathbf{0}$ is on the hyperplane through $\gamma(m)$ normal to
$-\gamma^{\prime}(m).$ Choose any point $p\in\partial B(\mathbf{0},1)\cap
O_{\gamma(m)}^{c}(-\gamma^{\prime}(m))$ and a shortest curve $\eta$ in
$O_{\gamma(m)}^{c}(-\gamma^{\prime}(m))=O^{c}$ from $\gamma(m)$ to $p$, in the
opposite direction of $\gamma$ at $\gamma(m).$ By Proposition 1, $\eta$ lies
in a 2-plane $X$ through $\gamma(m)$ and $p$, parallel to $\gamma^{\prime}(m)$
and it is a $C^{1,1}$ curve following a circular arc of length $\theta$ of
radius $1$ and a line segment to $p$. Let $A$ be the set $\{x\in\mathbf{R}%
^{n}:\left|  x\cdot\gamma^{\prime}(m)\right|  \leq1\}$ whose boundary consists
of two parallel hyperplanes. $\eta\subset A$, since $B(\mathbf{0},1)\cup
O_{\gamma(m)}(-\gamma^{\prime}(m))\subset int$ $A,$ $\gamma(m)$ and $p$ are in
$A.$ Consequently, $\eta\subset A\cap X\cap O^{c}.$ Since $\gamma(m)$ and $p$
are in different components of ($int$ $A)\cap X\cap O^{c},$ $\eta$ must pass
through $\partial A\cap X\cap O^{c}.$ This shows that $\theta\geq\frac{\pi}%
{2}$ and $\ell(\eta)>\frac{\pi}{2}.$ Suppose that $m\leq\frac{\pi}{2}.$ Then,
$\gamma([0,m])\subset O_{\gamma(m)}^{c}(-\gamma^{\prime}(m))$ by part
(a,$\frac{\pi}{2}$) and take $p=\gamma(0)\in\partial B(\mathbf{0},1).$ This
gives us a contradiction: $\frac{\pi}{2}<\ell(\eta)\leq\ell(\gamma
([0,m]))=m\leq\frac{\pi}{2}.$ Consequently, one must have $\frac{\pi}{2}<m,$
and $\frac{\pi}{2}<L-m$ by symmetry.

\textbf{(a:}$\pi$) By reparametrization and symmetry, it suffices to prove
this for $a=0$ and for $0\leq s\leq\pi.$ Suppose that $\gamma(b)\in
O_{\gamma(0)}(\gamma^{\prime}(0),1)$ for some $b\in(0,\pi]\cap I.$ Then
$\gamma(b)\in B(q,1)$ where $q=\gamma(0)+v$ for some unit vector $v$ normal to
$\gamma^{\prime}(0)$. There is a unique $c\in\lbrack0,b)$ such that
$\gamma((c,b])\subset B(q,1)$ and $\gamma(c)\in\partial B(q,1).$ One must have
$\gamma([0,b])\subset\overline{B(q,1)}$ by part (b), since $c<\pi,$
$\gamma(0)$ and $\gamma(c)$ are in $\partial B(q,1)$. $\gamma^{\prime}(c)$ is
tangent to $\partial B(q,1),$ since $\left\|  \gamma(t)-q\right\|  $ has a
local maximum at $t=c\neq0,$ and $c=0$ case is obvious. By part (a:$\frac{\pi
}{2}$), $\gamma$ must stay out of $O_{\gamma(c)}(\gamma^{\prime}(c),1)\supset
B(q,1)$ for $c\leq t\leq c+\frac{\pi}{2},$ which contradicts $\gamma
((c,b])\subset B(q,1).$ Hence, $\gamma([0,\pi]\cap\lbrack0,L])\cap
O_{\gamma(0)}(\gamma^{\prime}(0),1)=\emptyset.$

Assume that $\exists s_{0}\in I$ with $\gamma(s_{0})\in\partial O_{\gamma
(0)}^{c}(\gamma^{\prime}(0))$ and $0<s_{0}\leq\pi.$ Then, $\gamma(s_{0}) $ and
$\gamma(0)\in\partial B(q,1)$ where $q=\gamma(0)+v$ for some unit vector $v$
normal to $\gamma^{\prime}(0)$. By part (b) and the previous paragraph,
$(\gamma([0,s_{0}])\subset\overline{B(q,1)}\cap O_{\gamma(0)}^{c}%
(\gamma^{\prime}(0))$ which is the desired circle.

\textbf{(c)} Let $q=\gamma(a)+R\gamma^{\prime\prime}(a)$, and define
$f(s)=\frac{1}{2}\left\|  \gamma(s)-q\right\|  ^{2}.$

$f^{\prime}(s)=\gamma^{\prime}(s)\cdot\left(  \gamma(s)-q\right)  $ which is
lipschitz, and $f^{\prime}(a)=\gamma^{\prime}(a)\cdot\left(  -R\gamma
^{\prime\prime}(a)\right)  =0,$ by $\left\|  \gamma^{\prime}\right\|  \equiv1.$

$f^{\prime\prime}(s)=\gamma^{\prime\prime}(s)\cdot\left(  \gamma(s)-q\right)
+\gamma^{\prime}(s)\cdot\gamma^{\prime}(s)$ a.e., and $f^{\prime\prime
}(a)=\gamma^{\prime\prime}(a)\cdot\left(  -R\gamma^{\prime\prime}(a)\right)  +1<0.$

Hence, $\underset{s\rightarrow a^{+}}{lim}\frac{1}{s}(f^{\prime}(s)-f^{\prime
}(a))<0.$ There exists $\varepsilon>0$ such that $f^{\prime}(s)<0$ and
$f(s)<f(a)$, $\forall s\in(a,a+\varepsilon).$
\end{proof}
\end{proposition}

\begin{example}
$\pi$ in part (b) of the previous proposition is sharp. Consider the part of
the circle $(x-\varepsilon)^{2}+y^{2}=1$ outside the disc $x^{2}+y^{2}\leq1,$
in $\mathbf{R}^{2},$ for small $\varepsilon.$
\end{example}

\begin{lemma}
For all $C^{1,1}$ curves $\gamma:I\rightarrow\mathbf{R}^{n},$ analytic and
geometric focal distances are the same: $F_{g}(\gamma(I))=F_{k}(\gamma).$

\begin{proof}
Reparametrize $\gamma$ to assume that $\left\|  \gamma^{\prime}(s)\right\|
=1. $ $\frac{1}{F_{k}(\gamma)}\geq\kappa\gamma.$ By Proposition 2(a:$\frac
{\pi}{2}$) and rescaling, $\forall p\in\gamma,$ $\gamma$ locally avoids
$O_{p}(F_{k}(\gamma);\gamma)$ near $p$ and $F_{k}(\gamma)\leq F_{g}%
(p)=\inf\{r>0:p\in\overline{O_{p}(r;\gamma)\cap\gamma}\}.$ Hence,
$F_{k}(\gamma)\leq F_{g}(\gamma)=\inf_{p\in\gamma}F_{g}(p).$

Suppose that $F_{k}(\gamma)<F_{g}(\gamma),$ i.e. $\sup\kappa\gamma>\frac
{1}{F_{g}(\gamma)}.$ Define
\[
A=\left\{  s\in I:\kappa\gamma(s)>\frac{1}{F_{g}(\gamma)}\right\}  \text{ and
}B=\left\{  s\in I:\gamma^{\prime\prime}(s)\text{ exists}\right\}  .
\]
$A\neq\emptyset$ and the Lebesgue measure $\mu(B^{c})=0,$ where $X^{c}=I-X.$

\textbf{Case 1}. $A\cap B\neq\emptyset.$ There exists $s_{0}\in A\cap B$ such
that $c:=\left\|  \gamma^{\prime\prime}(s_{0})\right\|  =\kappa\gamma
(s_{0})>\frac{1}{F_{g}(\gamma)}.$ Choose $r$ such that $\frac{1}{c}%
<r<F_{g}(\gamma)$. Let $\eta(s)=c\gamma(\frac{s}{c}),$ so that $\left\|
\eta^{\prime}(s)\right\|  =1,\forall s,$ and $\left\|  \eta^{\prime\prime
}(cs_{0})\right\|  =1.$ By Proposition 2c, $\eta(cs_{0},cs_{0}+c\varepsilon
)\subset B(\eta(cs_{0})+cr\eta^{\prime\prime}(cs_{0}),cr)$ for some
$\varepsilon>0.$ Hence, $\gamma\left(  (s_{0},s_{0}+\varepsilon)\right)
\subset B(\gamma(s_{0})+r\frac{\gamma^{\prime\prime}(s_{0})}{\left\|
\gamma^{\prime\prime}(s_{0})\right\|  },r)\subset O_{\gamma(s_{0})}(r;\gamma).
$ However this contradicts $r<F_{g}(\gamma)$ by the definition of $F_{g}.$

\textbf{Case 2.} $A\cap B=\emptyset.$ Since $\gamma$ is $C^{1,1},$
$\gamma^{\prime}$ is absolutely continuous, $\gamma^{\prime\prime}(s)$ exists
almost everywhere by Rademacher's Theorem and $\left\|  \gamma^{\prime\prime
}(s)\right\|  =\kappa\gamma(s)\leq\frac{1}{F_{g}(\gamma)}$ a.e. By Lemma 1,
$\frac{1}{F_{k}(\gamma)}=\sup_{I}\kappa\gamma(s)\leq\frac{1}{F_{g}(\gamma)}$
which contradicts $F_{k}(\gamma)<F_{g}(\gamma).$

Neither of the cases is possible, hence one must have $F_{k}(\gamma
)=F_{g}(\gamma(I))$.
\end{proof}
\end{lemma}

\begin{definition}
Let $p,q\in\mathbf{R}^{n}$, $v\in UT\mathbf{R}_{p}^{n}$, $w\in UT\mathbf{R}%
_{q}^{n}$ and $\Lambda>0$ be given. Define $\mathcal{C}(p,q;v,w;\Lambda)$ to
be the set of all $C^{1,1}$ curves $\gamma:[0,L]\rightarrow\mathbf{R}^{n}$
with $\gamma(0)=p,$ $\gamma^{\prime}(0)=v,$ $\gamma(L)=q,$ $\gamma^{\prime
}(L)=w,$ $\left\|  \gamma^{\prime}\right\|  \equiv1,$ and $\kappa\gamma
\leq\Lambda,$ where $L=\ell(\gamma)$ is not fixed on $\mathcal{C}. $
\end{definition}

\begin{proposition}
There exists a shortest curve in $\mathcal{C}(p,q;v,w;\Lambda).$
\end{proposition}

\begin{proof}
Obviously, $\mathcal{C}(p,q;v,w;\Lambda)\neq\emptyset.$ Any sequence of curves
$\left\{  \gamma_{m}\right\}  _{m=1}^{\infty}$, with $\ell(\gamma
_{m})\rightarrow\inf\{\ell(\gamma):\gamma\in\mathcal{C}\}$ has uniformly
bounded lengths and all starting at $p$. Extend all $\gamma_{m}$ to a common
compact interval by following the lines $q+(s-\ell(\gamma_{m}))v$ after $q.$
By Lemma 1, $\forall\gamma\in\mathcal{C}$, $\left\|  \gamma^{\prime}%
(s)-\gamma^{\prime}(t)\right\|  \leq\Lambda\left|  s-t\right|  $, and thus,
$\mathcal{C}$\ is $C^{1}$-equicontinuous. $\mathcal{C}$\ is $C^{1}$-bounded by
$\left\|  \gamma^{\prime}\right\|  \equiv1$. $C^{0}$-equicontinuity and
boundedness are obvious. By Arzela-Ascoli Theorem, there exists a subsequence
of $\left\{  \gamma_{m}\right\}  _{m=1}^{\infty}$ uniformly converging to
$\gamma_{0}$ in $C^{1}$ sense: $(\gamma_{m}(s),\gamma_{m}^{\prime
}(s))\rightarrow$ $(\gamma_{0}(s),\gamma_{0}^{\prime}(s))$. $\gamma_{0}%
\in\mathcal{C}$, since all conditions of $\mathcal{C}$\ are preserved under
this convergence and $\ell(\gamma_{m})\rightarrow$ $\ell(\gamma_{0})$.
\end{proof}

\begin{proposition}
Let $\gamma:I=[0,L]\rightarrow\mathbf{R}^{n}$ be a shortest curve in
$\mathcal{C}(p,q;v,w;\Lambda).$ Then, $\forall s\in I,$ $(\kappa\gamma(s)=0 $
or $\Lambda).$ $\kappa\gamma^{-1}(\Lambda)$ is a closed subset of $I$, and
$\kappa\gamma^{-1}(0)$ is countable union of disjoint line segments.

\begin{proof}
By the upper semi-continuity of $\kappa\gamma,$ $\forall\lambda\leq\Lambda,$
$\kappa\gamma^{-1}([\lambda,\Lambda])$ is a closed subset of $I$ and
$J(\lambda)=\kappa\gamma^{-1}([0,\lambda))$ is countable union of relatively
open intervals in $I.$ Choose any $\lambda<\Lambda$ and $a<b$ in a given
component $J^{\prime}$ of $J(\lambda).$

Suppose that $\gamma^{\prime}(a)\neq\gamma^{\prime}(b).$ Choose any smooth
bump function $h:\mathbf{R}\rightarrow\lbrack0,1]$ such that \emph{supp}%
$(h)\subset\lbrack-1,1],$ $h(0)=1,$ and $\int\nolimits_{-1}^{1}h(s)ds=1.$ Let
$h_{n}$ be defined by $h_{n}(\frac{a+b}{2})=1$ and $h_{n}^{\prime
}(s)=n[h(n(s-a))-h(n(s-b))].$ Then,
\[
\lim_{n\rightarrow\infty}\int\nolimits_{I}h_{n}^{\prime}(s)\gamma^{\prime
}(s)ds=\gamma^{\prime}(b)-\gamma^{\prime}(a)\neq0.
\]
Choose and fix $n$ sufficiently large such that \emph{supp}$(h_{n})\subset
J^{\prime}$ and $-\int\nolimits_{J^{\prime}}h_{n}^{\prime}(s)\gamma^{\prime
}(s)ds=V\neq0.$ Let $\gamma_{\varepsilon}(s)=\gamma(s)+\varepsilon Vh_{n}(s)$
be a variation of $\gamma.$ By the First Variation formula, [CE, p6],
\[
\frac{d}{d\varepsilon}\ell(\gamma_{\varepsilon})|_{\varepsilon=0}%
=\int\nolimits_{I}[Vh_{n}(s)]^{\prime}\gamma^{\prime}(s)ds=-\left\|
V\right\|  ^{2}<0
\]
Hence, for sufficiently small $\varepsilon$, $\gamma_{\varepsilon}$ is
strictly shorter that $\gamma.$ For all $s<t:$
\begin{align*}
dil^{d}\gamma_{\varepsilon}^{\prime}(s,t)  &  =\frac{\left\|  \gamma
_{\varepsilon}^{\prime}(s)-\gamma_{\varepsilon}^{\prime}(t)\right\|  }%
{\ell_{st}(\gamma_{\varepsilon})}\leq\frac{\left\|  \gamma^{\prime}%
(s)-\gamma^{\prime}(t)\right\|  +\varepsilon\left\|  V\right\|  \left|
h_{n}^{\prime}(s)-h_{n}^{\prime}(t)\right|  }{t-s-\varepsilon\int
\nolimits_{s}^{t}\left|  h_{n}^{\prime}(u)\right|  \left\|  V\right\|  du}\\
&  \leq\frac{\left\|  \gamma^{\prime}(s)-\gamma^{\prime}(t)\right\|  }%
{t-s}+\varepsilon C(\left\|  V\right\|  ,\sup\left|  h_{n}^{\prime}\right|
,\sup\left|  h_{n}^{\prime\prime}\right|  )
\end{align*}
By Remark 1.i, for sufficiently small $\varepsilon$, $\kappa\gamma
_{\varepsilon}\leq\frac{\Lambda+\lambda}{2}<\Lambda,$ and $\gamma
_{\varepsilon}\in\mathcal{C}.$ This contradicts the minimality of $\gamma. $
Consequently, $\gamma^{\prime}$ is constant on $J^{\prime}.$ $\forall
\lambda<\Lambda,$ $\gamma(J(\lambda))$ is a countable union of disjoint line
segments, to conclude that $\gamma(J(\Lambda))$ is a countable union of
disjoint line segments, and $\kappa\gamma(J(\Lambda))\equiv0.$
\end{proof}
\end{proposition}

\subsection{Proof of Theorem 1}

\begin{proof}
By using dilations of $\mathbf{R}^{n}$, one can assume that $\Lambda=1.$ We
are going to proceed in proving parts (a) and (b) simultaneously, and point
out the differences when they are needed. Let $A=\pi$ for part (a) and
$A=2\pi$ for part (b). By Proposition 4, there exist maximally chosen $c$ and
$d$ such that $s_{0}\in\lbrack c,d]\subset\lbrack0,L]$ and $\gamma([c,d])$ is
a line segment $L_{0}$.

Assume that $\gamma([a,b])$ is a $CLC(1)$-curve for $[c,d]\subset\lbrack
a,b]\subset\lbrack0,L]$. We will show that if $a>0$ and $c-a<A$, then
$\exists\delta>0$ such that $\gamma([a-\delta,b])$ still is a $CLC(1)$ curve.

For $r\in\lbrack0,a]$ define $J_{r}=\{-\lambda\gamma^{\prime}(a-r):\lambda
>0\}$ and $V_{r}=O_{\gamma(a-r)}^{c}(\gamma^{\prime}(a-r)). $ Choose
$d^{\prime}=d$ when $\gamma((c,d])\cap J_{0}=\emptyset;$ otherwise, by
$\gamma(d^{\prime})\in\gamma((c,d])\cap J_{0}\neq\emptyset.$ Let
$\varepsilon=\frac{1}{2}\min(d^{\prime}-c,A-(c-a)),$ $c_{1}=c+\varepsilon,$
and $m=\gamma(c_{1}).$ $\ m\in intV_{0}$, since $\gamma([a,c])$ is an arc of a
circle of radius 1 and $\gamma([c,c_{1}])$ is a line segment, $a\leq c<c_{1},
$ and Proposition 2a. One obtains that $\forall r\in\lbrack0,\varepsilon),$
$\gamma([a-r,c_{1}])\subset V_{r}$ by Proposition 2a and $c_{1}%
-(a-r)<c-a+2\varepsilon\leq\pi$ for part (a), and by $R_{O}(\gamma)\geq1$ for
part (b). $m\in intV_{r}$, since $\gamma([c,c_{1}])=L_{0}$ is a line segment.

For each fixed $r\in\lbrack0,\varepsilon),$ define $\gamma_{r}$ to be a
shortest curve parametrized by arclength \ from $\gamma(a-r)$ to $m$ within
$V_{r}$ \textbf{without curvature restrictions, }by using Proposition 1.
$\gamma_{r}$ follows a circular arc of radius $1$ starting from $\gamma(a-r) $
along $\partial V_{r}$, then a line segment $L_{r}$ of positive length until
$m.$ In Proposition 1, $\gamma_{r}^{\prime}$ at $m$ is not controlled.

\textbf{Claim 1. }$\gamma_{r}^{\prime}(m)=\gamma^{\prime}(m)$ for sufficiently
small $r>0.$

$\exists\delta_{1}>0$ such that $\forall r\in\lbrack0,\delta_{1}),$
$\ell(L_{r})\geq\frac{\varepsilon}{2}$, $d(m,J_{r})>0,$ since $J_{r}$ and
$\partial V_{r}$ change continuously in $r$, $\ell(L_{0})=\varepsilon$ and
$d(m,J_{0})>0.$ For $r\in\lbrack0,\delta_{1}),$ $\gamma_{r}$ is uniquely
defined. $\lim_{r\rightarrow0^{+}}\measuredangle_{m}(L_{0},L_{r})=0,$
otherwise one can construct a shortest curve other than $\gamma$ from
$\gamma(a)$ to $m$ in $V_{0}$ contradicting Proposition 1b$.$ $\exists
\delta_{2}>0$ such that $\forall r\in\lbrack0,\delta_{2}),$ $\measuredangle
_{m}(L_{0},L_{r})\leq2\tan^{-1}\frac{\varepsilon}{2}.$

Let $\delta=\min(\varepsilon,\delta_{1},\delta_{2}).$ $\forall r\in
\lbrack0,\delta),$ define a curve $\tilde{\gamma}_{r}$ which follows $\gamma$
from $p$ to $\gamma(a-r),$ then $\gamma_{r}$ from $\gamma(a-r)$ to $m$, and
$\gamma$ from $m$ to $q.$ $\tilde{\gamma}_{r}$ is $C^{1}$ at $\gamma(a-r)$,
squeezed by $O_{\gamma(a-r)}(\gamma^{\prime}(a-r)).$ Recall that
$\varepsilon\leq\frac{d-c}{2}$ and $L_{0}=\gamma([c,c_{1}]).$ Define the line
segment $L_{0}^{\prime}:=\gamma([c_{1},c_{1}+\frac{\varepsilon}{2}]).$

Suppose that $\tilde{\gamma}_{r}$ is not $C^{1}$ at $m$ for some
$r\in(0,\delta),$ that is $\measuredangle_{m}(L_{0},L_{r})=\pi-\measuredangle
_{m}(L_{0}^{\prime},L_{r}):=\alpha>0.$ Fix such an $r.$ $\gamma\cap
B(m,\frac{\varepsilon}{2})$ is a union of two segments of length
$\frac{\varepsilon}{2}$, joined at $m$ with an angle of $\pi-\alpha$, in
$L_{0}^{\prime}\cup L_{r}.$ There exists a unique circle $C$ of radius $1$ in
the same 2-plane as $L_{0}^{\prime}\cup L_{r}$, tangent to $L_{r}$ at $p_{1}$
and tangent to $L_{0}^{\prime}$ at $p_{2}$ where $\left\|  p_{i}-m\right\|
\leq\frac{\varepsilon}{2},$ since $\alpha\leq2\tan^{-1}\frac{\varepsilon}{2}.
$ Let $\tilde{\gamma}$ be the $C^{1}$ curve obtained from $\tilde{\gamma}_{r}
$ by replacing $L_{0}^{\prime}\cup L_{r}$ between $p_{1}$ and $p_{2}$ by the
shorter arc of $C$ between $p_{1}$ and $p_{2}.$%
\[
\ell(\tilde{\gamma})<\ell(\tilde{\gamma}_{r})\leq\ell(\gamma)\text{ and
}\kappa\tilde{\gamma}\leq1
\]
This contradicts the minimality of $\gamma$ in $\mathcal{C}$. Hence, $\forall
r\in\lbrack0,\delta),$ $\tilde{\gamma}_{r}$ is $C^{1}$ at $m,$ and
$\measuredangle_{m}(L_{0},L_{r})=0.$ This proves Claim 1.

For each given $r\in(0,\delta):$

1. $\tilde{\gamma}_{r}\in C^{1}$ and $\kappa\tilde{\gamma}_{r}\leq1,$ hence
$\tilde{\gamma}_{r}\in\mathcal{C}$ and $\ell(\tilde{\gamma}_{r})\geq
\ell(\gamma).$

2. $\gamma$ and $\tilde{\gamma}_{r}$ follow the same path before $\gamma(a-r)
$ as well as after $m.$

3. $\gamma([a-r,c_{1}])\subset V_{r},$ $\gamma_{r}$ is the unique shortest
curve from $\gamma(a-r)$ to $m$ in $V_{r}$, and hence $\ell(\gamma
([a-r,c_{1}]))\geq\ell(\gamma_{r}).$

Consequently, $\ell(\tilde{\gamma}_{r})=\ell(\gamma),$ $\gamma$ and
$\tilde{\gamma}_{r}$ are equal up to parametrization, and $\forall r\in
\lbrack0,\delta),$ $\gamma|[a-r,b])$ is a $CLC(1)$-curve. Obviously, this
extends to $[a-\delta,b]$ and to $[a-\delta,b+\delta^{\prime}]$ for some
$\delta^{\prime}>0$, by symmetry \ when $b<L$ and $d-b<A.$ One chooses $a_{0}$
and $b_{0}$ maximally so that $[c,d]\subset\lbrack a_{0},b_{0}]\subset
\lbrack0,L]$ and $\gamma|[a_{0},b_{0}]$ is a $CLC(1)$ curve. It follows from
the construction of $\delta$ that:

Part (a): $(0=a_{0}$ or $c-a_{0}\geq A=\pi)$ and $(L=b_{0}$ or $d-b_{0}\geq
A=\pi)$

Part (b): $0=a_{0}$ and $L=b_{0},$ since $c-a_{0}=A=2\pi$ case creates a
complete circle through $\gamma(c),$ which contradicts the minimality of
$\gamma.$
\end{proof}

\section{Relatively Extremal Knots in $\mathbf{R}^{n}$}

A knot class $[\theta]$ is a free $C^{0}-$homotopy class of embeddings of
$\gamma:\mathbf{S}^{1}\rightarrow\mathbf{R}^{n}.$ In this section$,$
$\gamma:\mathbf{S}^{1}\rightarrow\mathbf{R}^{n}$ denotes a simple$-C^{1}%
-$closed curve, by identifying $\mathbf{S}^{1}\cong\mathbf{R}/L\mathbf{Z}$ and
$K=image(\gamma)$. In other words, $\gamma(t+L)=\gamma(t)$ and $\gamma
^{\prime}(t+L)=\gamma^{\prime}(t)$, $\forall t\in\mathbf{R}$ with $\left\|
\gamma^{\prime}\right\|  \neq0$ and $\gamma$ is one-to-one on $[0,L).$
Interval notation will be used to describe subsets of $\mathbf{R}/L\mathbf{Z.}$

\begin{definition}
For any simple$-C^{1,1}-$closed curve $\gamma:\mathbf{S}^{1}\rightarrow
\mathbf{R}^{n}$, one defines the ropelength or extrinsically isoembolic length
to be $\ell_{e}(\gamma)=\frac{\ell(\gamma)}{R_{o}(\gamma)}=\frac
{vol_{1}(\gamma(S^{1}))}{i(\gamma(S^{1}),\mathbf{R}^{n})}.$
\end{definition}

\begin{definition}
i. A simple$-C^{1,1}-$closed curve $\gamma_{0}$ is called an ideal (thickest)
knot in $[\theta],$ if $\ell_{e}(\gamma_{0})\leq\ell_{e}(\gamma),\forall
\gamma\in\lbrack\theta]\cap C^{1,1}.$

ii. $\gamma_{0}$ is called relatively extremal, if there exists an open set
$\mathcal{U}$ in $C^{1}$ topology such that $\gamma_{0}\in\mathcal{U}$ and
$\ell_{e}(\gamma_{0})\leq\ell_{e}(\gamma),\forall\gamma\in\mathcal{U}%
\cap\lbrack\theta]\cap C^{1,1}.$
\end{definition}

We consider two curves $\gamma_{1}$ and $\gamma_{2}$ to be geometrically
equivalent if there exists an orientation preserving $h:\mathbf{R}%
^{n}\rightarrow\mathbf{R}^{n}$, a composition of an isometry and a dilation
$(\mathbf{x}\rightarrow\lambda\mathbf{x,}$ $\lambda\neq0)$ of $\mathbf{R}^{n}%
$, such that $h(\gamma_{1})=\gamma_{2}$ up to a reparametrization. On each
geometric equivalence class of $C^{1,1}-$closed curves, $l_{e}$ remains constant.

\begin{theorem}
(Thickness Formula) For every simple$-C^{1,1}-$closed curve $\gamma$ in
$\mathbf{R}^{n},$ and $K=image(\gamma)$, one has $i(K,M)=R_{O}(K,M)=\min
\{F_{k}(\gamma),\frac{1}{2}MDC(K)\}.$
\end{theorem}

\begin{proof}
See [L], for $n=3$ case. This is a consequence of Thickness Formula [D,
\emph{Theorem 1}] and Lemma 2.
\end{proof}

\begin{proposition}
Let $\left\{  \gamma_{m}\right\}  _{m=1}^{\infty}:\mathbf{S}^{1}\left(
\cong\mathbf{R}/L\mathbf{Z}\right)  \rightarrow\mathbf{R}^{n}$ be a sequence
uniformly converging to $\gamma_{0}$ in $C^{1}$ sense, i.e. $(\gamma
_{m}(s),\gamma_{m}^{\prime}(s))\rightarrow$ $(\gamma_{0}(s),\gamma_{0}%
^{\prime}(s))$ uniformly on $\mathbf{S}^{1}.$ Let $K_{m}=\gamma_{m}%
(\mathbf{S}^{1}).$

i. If $R_{O}(K_{m})\geq r$ for sufficiently large $m$, then $R_{O}(K_{0})\geq
r$. Consequently, $\lim\sup_{m}R_{O}(K_{m})\leq R_{O}(K_{0}).$

ii. If $\lim\inf_{m}MDC(K_{m})>0,$ then $\lim\inf_{m}MDC(K_{m})\geq MDC(K_{0}).$
\end{proposition}

\begin{proof}
i. Suppose that $R_{O}(K_{0})<r$, for a given $r>0.$ By the definition
$R_{O},$ there exists $a\in\mathbf{S}^{1}$, $v\in\mathbf{R}^{n}$ with
$\left\|  v\right\|  =1$ and $v\cdot\gamma_{0}^{\prime}(a)=0$ such that
$B(\gamma_{0}(a)+rv,r)\cap K_{0}\neq\emptyset.$ One can find $\gamma_{0}(b)\in
B(\gamma_{0}(a)+rv,r-\varepsilon)$ for sufficiently small $\varepsilon>0.$
Choose a sequence $\left\{  v_{m}\right\}  _{m=1}^{\infty}$ in $\mathbf{R}%
^{n}$ such that $\forall m,\left\|  v_{m}\right\|  =1$, $v_{m}\cdot\gamma
_{m}^{\prime}(a)=0,$ and $v_{m}\rightarrow v$. Then for sufficiently large
$m$,
\[
\left\|  \left(  \gamma_{m}(a)+rv_{m}\right)  -\gamma_{m}(b)\right\|
<\left\|  \left(  \gamma_{0}(a)+rv\right)  -\gamma_{0}(b)\right\|
+\varepsilon<r.
\]
Hence, $B(\gamma_{m}(a)+rv_{m},r)\cap\gamma_{m}\neq\emptyset$ and $R_{O}%
(K_{m})<r,$ for sufficiently large $m,$ which contradicts the hypothesis.
Consequently, $R_{O}(K_{0})\geq r.$

ii. We will use the same indices for subsequences. Let $a=\lim\inf
_{m}MDC(K_{m}),$ and choose a subsequence with $a=\lim_{m}MDC(K_{m})$ and
$MDC(K_{m})>0,\forall m.$ By compactness of $K_{m}$ and positivity of
$MDC(K_{m})$, there exists a minimal double critical pair $\{p_{m},q_{m}\}$
for $K_{m}$, $\ell(\gamma_{p_{m}q_{m}})=MDC(K_{m}),\forall m.$ Since $K_{0} $
is compact and $a>0,$ there exists subsequences $p_{m}\rightarrow p_{0}\in
K_{0}$, $q_{m}\rightarrow q_{0}\in K_{0},$ and $\gamma_{p_{m}q_{m}}%
\rightarrow\gamma_{p_{0}q_{0}}$ in $C^{1}$ sense. Line segments converge to
line segments, and normality to $C^{1}$ curves is preserved under $C^{1}$
limits. Hence $\{p_{m},q_{m}\}$ is a double critical pair for $K_{0}.$%
\[
MDC(K_{0})\leq\ell(\gamma_{p_{0}q_{0}})=\lim_{m}\ell(\gamma_{p_{m}q_{m}}%
)=\lim_{m}MDC(K_{m})=a.
\]
\end{proof}

\begin{definition}
Let $\gamma:\mathbf{S}^{1}\rightarrow\mathbf{R}^{n}$ be a simple$-C^{1,1}%
-$closed curve in $\mathbf{R}^{n},$ with $\gamma^{\prime}\neq0.$ Define

i. $I_{c}=\{x\in\mathbf{S}^{1}:\exists y\in\mathbf{S}^{1}$ such that $\left\|
\gamma(x)-\gamma(y)\right\|  =MDC(\gamma)$ and

$\qquad\left(  \gamma(x)-\gamma(y)\right)  \cdot\gamma^{\prime}(x)=\left(
\gamma(x)-\gamma(y)\right)  \cdot\gamma^{\prime}(y)=0\}$ and $K_{c}=\gamma
_{c}=\gamma(I_{c})$

ii. $I_{z}=\{x\in\mathbf{S}^{1}:\kappa\gamma(x)=0\}$ and $K_{z}=\gamma
_{z}=\gamma(I_{z})$

iii. $I_{mx}=\{x\in\mathbf{S}^{1}:\kappa\gamma(x)=1/R_{O}(\gamma)\}$ and
$K_{mx}=\gamma_{mx}=\gamma(I_{mx})$

iv. $I_{b}=\{x\in\mathbf{S}^{1}:0<\kappa\gamma(x)<1/R_{O}(\gamma)\}$ and
$K_{b}=\gamma_{b}=\gamma(I_{b})$
\end{definition}

\begin{remark}
$K_{c}$ and $K_{mx}$ are closed subsets of $K$. This is obvious for $K_{c}$ by
the continuity of $\gamma^{\prime}.$ See the proof of Proposition 8, for
$K_{mx}.$
\end{remark}

\begin{proposition}
For any knot class $[\theta]$ in $\mathbf{R}^{n}$, $\exists\gamma_{0}%
\in\lbrack\theta]\cap C^{1,1}$ such that

i. $\forall\gamma\in\lbrack\theta]\cap C^{1,1},$ $0<\ell_{e}(\gamma_{0}%
)\leq\ell_{e}(\gamma),$ and hence

ii. $\forall\gamma\in\lbrack\theta]\cap C^{1,1},\left(  \ell(\gamma_{0}%
)=\ell(\gamma)\Longrightarrow R_{O}(\gamma_{0})\geq R_{O}(\gamma)\right)  .$
\end{proposition}

\begin{proof}
Let $\mathbf{T}_{x}=\{\gamma\in\lbrack\theta]\cap C^{1,1}:\gamma
(0)=\mathbf{0},$ $\left\|  \gamma^{\prime}\right\|  \equiv1,$ $\ell(\gamma)=1$
and $R_{O}(\gamma)\geq x\}.$ Every geometric equivalence class of $C^{1,1}%
-$closed curves has a representative in $\mathbf{T}_{0}$. $\forall\gamma
\in\mathbf{T}_{0}$, $MDC(\gamma)\leq\frac{1}{2}.$ Choose any $\gamma_{1}%
\in\lbrack\theta]\cap C^{\infty}\cap\mathbf{T}_{0}$ and set $R_{O}(\gamma_{1})=A.$

By the Thickness Formula, $0<A\leq M:=\sup\{R_{O}(\gamma):\gamma\in
\mathbf{T}_{A}\}\leq\frac{1}{4}$. $\forall\gamma\in\mathbf{T}_{A},F_{k}%
(\gamma)\geq A$ and hence $\kappa(\gamma)\leq\frac{1}{A},$ or equivalently,
$\left\|  \gamma^{\prime}(s)-\gamma^{\prime}(t)\right\|  \leq\frac{\left|
s-t\right|  }{A},\forall s,t\in\mathbf{R/Z}$. Consequently, $\mathbf{T}_{A}%
$\ is $C^{1}$-equicontinuous and $C^{1}$-bounded: $\left\|  \gamma^{\prime
}\right\|  \equiv1$. $C^{0}$-equicontinuity and boundedness is obvious.

For all $k\in\mathbf{N}^{+}$, there exists $\gamma_{k}\in\mathbf{T}_{A}$ such
that $R_{O}(\gamma_{k})\geq M-\frac{1}{k}.$ By Arzela-Ascoli Theorem, there
exists a subsequence of $\left\{  \gamma_{m}\right\}  _{m=1}^{\infty}$
uniformly converging to $\gamma_{0}$ in $C^{1}$ sense: $(\gamma_{m}%
(s),\gamma_{m}^{\prime}(s))\rightarrow$ $(\gamma_{0}(s),\gamma_{0}^{\prime
}(s))$. By Proposition 5, $R_{O}(\gamma_{0})\geq M.$ The rest is
straightforward: $\gamma_{0}\in\mathbf{T}_{M}\subset\mathbf{T}_{A}%
\subset\mathbf{T}_{0},$ $R_{O}(\gamma_{0})=M=\sup\{R_{O}(\gamma):\gamma
\in\mathbf{T}_{A}\}=\sup\{R_{O}(\gamma):\gamma\in\mathbf{T}_{0}\}$, and
$\ell_{e}(\gamma_{0})=\frac{1}{M}=\inf\{\ell_{e}(\gamma):\gamma\in
\mathbf{T}_{0}\}.$ Finally, all curves in the geometric equivalence class of
$\gamma_{0}$ are $\ell_{e}$-minimizers in $[\theta]\cap C^{1,1}.$
\end{proof}

\begin{proposition}
Let $\left\{  \gamma_{m}\right\}  _{m=1}^{\infty}:\mathbf{S}^{1}%
\rightarrow\mathbf{R}^{n}$ be a sequence uniformly converging to $\gamma$ in
$C^{1}$ sense, $K=\gamma(\mathbf{S}^{1})$ and $K_{m}=\gamma_{m}(\mathbf{S}%
^{1})$ satisfying

a. $\exists C<\infty,\forall m,$ $\sup\kappa\gamma_{m}\leq C,$ and

b. $\exists$ compact $A\subset\mathbf{S}^{1}$ such that $\{s\in\mathbf{S}%
^{1}:\gamma_{m}(s)\neq\gamma(s)\}\subset A,\forall m$.

Then both of the following hold.

i. If $A\cap I_{c}=\emptyset,$ then $\exists m_{1}\forall m\geq m_{1}%
,(MDC(K_{m})\geq MDC(K)).$

ii. If $F_{k}(\gamma)<\frac{1}{2}MDC(K)$ and $(F_{k}(\gamma_{m})\geq
F_{k}(\gamma)),\forall m$,

\qquad then $\exists m_{1}\forall m\geq m_{1},(R_{O}(K_{m})\geq R_{O}(K)).$
\end{proposition}

\begin{proof}
All subsequences will be denoted by the same index $m$. The critical pairs
will be identified from the domain $\mathbf{S}^{1}.$

i. Suppose there exists a subsequence $\gamma_{m}$ such that $\forall
m(MDC(K_{m})<MDC(K))$. For all $m$, there exists a minimal double critical
pair $\left\{  x_{m},y_{m}\right\}  $ in $\mathbf{S}^{1}$ for $\gamma_{m},$
$MDC(K_{m})=\left|  \gamma_{m}(x_{m})-\gamma_{m}(y_{m})\right|  <MDC(K).$
Then, $\forall m(d_{\mathbf{S}^{1}}(x_{m},y_{m})\geq\frac{\pi}{C})$ by
Proposition 2b. There exist subsequences $x_{m}\rightarrow x_{0}$,
$y_{m}\rightarrow y_{0}$ and $d_{\mathbf{S}^{1}}(x_{0},y_{0})\geq\frac{\pi}%
{C}.$ By the uniform convergence of $\gamma_{m}^{\prime}\rightarrow
\gamma^{\prime},$ $\left\{  x_{0},y_{0}\right\}  $ is a double critical pair
for $\gamma.$
\[
MDC(K)\leq\left|  \gamma(x_{0})-\gamma(y_{0})\right|  =\lim_{m}\left|
\gamma_{m}(x_{m})-\gamma_{m}(y_{m})\right|  =\lim_{m}MDC(K_{m})\leq MDC(K)
\]
Hence, $\left\{  x_{0},y_{0}\right\}  $ is a minimal double critical pair for
$\gamma$ and $\left\{  x_{0},y_{0}\right\}  $ $\subset I_{c}.$ Since $I_{c}$
and $A$ are disjoint compact subsets of $\mathbf{S}^{1}$, the subsequences
$\left\{  x_{m}\right\}  _{m=1}^{\infty}$ and $\left\{  y_{m}\right\}
_{m=1}^{\infty}$ can be taken in $\mathbf{S}^{1}-A.$ $\forall m,$ $\left\{
x_{m},y_{m}\right\}  $ is a double critical pair for $\gamma,$ since
$\gamma_{m}=\gamma$ on $\mathbf{S}^{1}-A.$%
\[
MDC(K)\leq\left|  \gamma(x_{m})-\gamma(y_{m})\right|  =\left|  \gamma
_{m}(x_{m})-\gamma_{m}(y_{m})\right|  =MDC(K_{m})
\]
which contradicts the initial assumption. Consequently, there does not exist
any subsequence $\gamma_{m}$ such that $\forall m,(MDC(K_{m})<MDC(K))$,
proving (i).

ii. $MDC(K)>2F_{k}(\gamma)=2R_{O}(K)$ and $\forall m$, $(F_{k}(\gamma_{m})\geq
F_{k}(\gamma))$. Suppose that there exists a subsequence $\gamma_{m}$ such
that $\forall m(MDC(K_{m}))<2R_{O}(K).$ For all $m$, there exists a minimal
double critical pair $\left\{  x_{m},y_{m}\right\}  $ in $\mathbf{S}^{1}$ for
$\gamma_{m},$ $MDC(K_{m})=\left|  \gamma_{m}(x_{m})-\gamma_{m}(y_{m})\right|
<2R_{O}(K).$ Then as in part (i), $\forall m(d_{\mathbf{S}^{1}}(x_{m}%
,y_{m})\geq\frac{\pi}{C})$ and by taking subsequences $x_{m}\rightarrow x_{0}%
$, $y_{m}\rightarrow y_{0}$, $d_{\mathbf{S}^{1}}(x_{0},y_{0})\geq\frac{\pi}%
{C},$ one obtains a double critical pair $\left\{  x_{0},y_{0}\right\}  $ for
$\gamma.$%
\[
MDC(K)\leq\left|  \gamma(x_{0})-\gamma(y_{0})\right|  =\lim_{m}\left|
\gamma_{m}(x_{m})-\gamma_{m}(y_{m})\right|  =\lim_{m}MDC(K_{m})\leq2R_{O}(K)
\]
which contradicts the hypothesis. Hence, $\exists m_{1}\forall m\geq
m_{1},(MDC(K_{m})\geq2R_{O}(K)),$ to conclude that
\[
R_{O}(K_{m})=\min\left(  F_{k}(\gamma_{m}),\frac{1}{2}MDC(K_{m})\right)
\geq\min\left(  F_{k}(\gamma),R_{O}(K)\right)  =R_{O}(K).
\]
\end{proof}

\begin{proposition}
(Also see [GM, p11, 12] for another version for smooth ideal knots.) Let
$\gamma$ be a relatively extremal knot.

i. If $MDC(K)=2R_{O}(K),$ then $K-(K_{c}\cup K_{mx})$ is a countable union of
open ended line segments, and hence $I_{b}\subset I_{c}.$

ii. If $MDC(K)>2R_{O}(K)$, then $K-K_{mx}$ is a countable union of open ended
line segments.
\end{proposition}

\begin{remark}
Theorem 2 shows that $K-K_{mx}$ is actually empty when $MDC(K)>2R_{O}(K).$
\end{remark}

\begin{proof}
Let $\mathcal{U}$ be an open set in $C^{1}$ topology such that $\gamma
\in\mathcal{U}$ and $\ell_{e}(\gamma)\leq\ell_{e}(\eta),\forall\eta
\in\mathcal{U}\cap\lbrack\theta]\cap C^{1,1}.$

i. Let $\Lambda=\sup\kappa\gamma.$ As in the proof of Proposition 4, for all
$\lambda\leq\Lambda,$ $\kappa\gamma^{-1}([0,\lambda))-I_{c}$ is countable
union of relatively open intervals in $\mathbf{S}^{1}(=\mathbf{R}%
/L\mathbf{Z)}.$ Choose any $\lambda<\Lambda$ and a closed interval $\left[
a,b\right]  $ contained in a component of $\kappa\gamma^{-1}([0,\lambda
))-I_{c}.$ By repeating the proof of Proposition 4, if $\gamma|\left[
a,b\right]  $ is not a line segment, then there exists a length decreasing
variation $\gamma_{\varepsilon}(s)=\gamma(s)+\varepsilon Vh_{n}(s)$ supported
in $\left[  a,b\right]  $. There exists a sufficiently small $\varepsilon
_{1}>0$ such that $\forall\varepsilon,0<\varepsilon\leq\varepsilon_{1},$ one has

\qquad1. $\gamma_{\varepsilon}$ and $\gamma$ belong to the same knot class and
$\gamma_{\varepsilon}\in\mathcal{U.}$

\qquad2. $\ell(\gamma_{\varepsilon})<\ell(\gamma),$ (proof of Proposition 4)

\qquad3. $\kappa\gamma_{\varepsilon}\leq\Lambda$ and hence $F_{k}%
(\gamma_{\varepsilon})\geq F_{k}(\gamma),$ (proof of Proposition 4), and

\qquad4. $MDC(K_{\varepsilon})\geq MDC(K)$ (Proposition 7(i)\ and $\left[
a,b\right]  \cap I_{c}=\emptyset$).

By the Thickness Formula, one obtains $R_{O}(K_{\varepsilon})\geq R_{O}(K)$
and $\ell_{e}(\gamma_{\varepsilon})=\frac{\ell(\gamma_{\varepsilon})}%
{R_{O}(K_{\varepsilon})}<\frac{\ell(\gamma)}{R_{O}(K)}=\ell_{e}(\gamma)$ which
is in contradiction with the hypothesis. Hence, $\gamma|\left[  a,b\right]  $
must be a line segment. Consequently, $I_{b}-I_{c}=\emptyset.$

ii. $MDC(K)>2R_{O}(K)=2F_{k}(\gamma).$ The proof is essentially the same as in
(i), with the following modifications. $\left[  a,b\right]  $ is taken in any
component of $\kappa\gamma^{-1}([0,\lambda)),$ thus $\left[  a,b\right]  \cap
I_{c}(\gamma)$ may not be empty. 1-3 above hold. To conclude $R_{O}%
(K_{\varepsilon})\geq R_{O}(K)$, one uses Proposition 7(ii). In this case,
$I_{b}=\emptyset$ and $K-K_{mx}$ is a countable union of open ended line segments.
\end{proof}

\subsection{Proof of Theorem 2}

\begin{proof}
Let $\mathcal{U}$ be an open set in $C^{1}$ topology such that $\gamma
\in\mathcal{U}$ and $\ell_{e}(\gamma)\leq\ell_{e}(\eta),\forall\eta
\in\mathcal{U}\cap\lbrack\theta]\cap C^{1,1}.$ We prove part (ii) first.

By Proposition 8, there exist maximally chosen $a,b$ such that $\gamma|(a,b) $
is an open ended line segment, $s_{0}\in(a,b)$ and $\left(  a,b\right)  \cap
I_{c}=\emptyset.$ If $b\in I_{c},$ then take $d=b,$ to finish the positive
direction. If $b\notin I_{c}$, proceed as follows. Assume that $\gamma
|[s_{0},b+\varepsilon]$ is a $CLC(F_{k}(\gamma)^{-1})$-curve (in fact, line
segment followed by circular arc) such that $0\leq\varepsilon<\pi F_{k}%
(\gamma)$ and $[s_{0},b+\varepsilon]\cap I_{c}=\emptyset.$ We will show that
the same is true for some $\varepsilon_{1}>\varepsilon.$ We point out that a
priori $\gamma|[s_{0},b+\varepsilon_{1}]$ is not known to be a shortest curve
in a certain $\mathcal{C,}$ replacing it with a shortest curve may create a
knot outside $\mathcal{U}$ or the knot class of $\gamma$, and this shortest
curve may not have a point of zero curvature.

Let $\left\{  b_{m}\right\}  _{m=1}^{\infty}$ be a sequence and $A>0$ be such that

\qquad1. $b_{m+1}<b_{m},\forall m\in\mathbf{N}^{+},$

\qquad2. $b_{m}\rightarrow b_{0}=b+\varepsilon,$

\qquad3. $\varepsilon<b_{m}-b\leq A\leq\pi F_{k}(\gamma),$ and

\qquad4. $[s_{0},b+A]\cap I_{c}=\emptyset.$

Since $\gamma(s_{0})\in int$ $O_{\gamma(b_{0})}^{c}(-\gamma^{\prime}(b_{0})),
$ $\gamma(s_{0})\in int$ $O_{\gamma(b_{m})}^{c}(-\gamma^{\prime}(b_{m}))$ for
sufficiently large $m\geq m_{0}.$ \ Let $f_{m}(s)$ be the unique shortest
curve parametrized by arclength in $\mathcal{C}(\gamma(s_{0}),\gamma
(b_{m});\gamma^{\prime}(s_{0}),\gamma^{\prime}(b_{m});F_{k}(\gamma)^{-1})$ by
Proposition 3, such that $f_{m}(s_{0})=\gamma(s_{0})$ and $f_{m}(c_{m}%
)=\gamma(b_{m}).$ Extend $f_{m} $ to $[s_{0},b+A]$ in a $C^{1}$ fashion beyond
$\gamma(b_{m})$ by $\gamma(s-c_{m}+b_{m})=$ $f_{m}(s).$

$\forall m,$ $\kappa f_{m}\leq F_{k}(\gamma)^{-1},\left\|  f_{m}^{\prime
}\right\|  =1$ and $f_{m}(s_{0})=\gamma(s_{0}).$ Hence, the sequence $\left\{
f_{m}\right\}  _{m=1}^{\infty}$ is $C^{1}$ equicontinuous and bounded. By
Arzela-Ascoli Theorem, there exists a convergent subsequence (which we will
denote by the same subindices $m)$ $f_{m}\rightarrow f_{0}$ uniformly in
$C^{1}$ topology. By the construction above, $f_{0}$ follows $\gamma$ past
$\gamma(b_{0})$ and $f_{0}(c_{0})=\gamma(b_{0})$ for some $c_{0.}$%
\begin{align*}
c_{m}-s_{0}  &  \leq b_{m}-s_{0}\\
\underset{m}{\lim\sup}\text{ }c_{m}  &  \leq b_{0}\\
c_{0}  &  \leq b_{0}%
\end{align*}%
\begin{align*}
f_{0}(s_{0})  &  =\gamma(s_{0})\text{ and }f_{0}(c_{0})=\gamma(b_{0})\\
f_{0}  &  \in C^{1}\text{ and }f_{0}^{\prime}(c_{0})=\gamma^{\prime}(b_{0})\\
\kappa f_{0}  &  \leq F_{k}(\gamma)^{-1}%
\end{align*}
By Proposition 1, $\gamma$ which is a line segment followed by a circular arc
is the unique shortest curve satisfying the last 3 conditions. Consequently,
$b_{0}=c_{0}$ and $f_{0}=\gamma$ on $[s_{0},b+A].$

Let $\gamma_{m}$ be the curve obtained from $\gamma$ by replacing
$\gamma|\left[  s_{0},b_{m}\right]  $ by $f_{m}|[s_{0},c_{m}].$ Reparametrize
$\gamma_{m}$ (not necessarily with respect to arclength) so that $\gamma
_{m}(s)=\gamma(s)$ for $s\notin\lbrack s_{0},b+A]$ and $\gamma_{m}%
\rightarrow\gamma$ in $C^{1}$ sense on $\mathbf{S}^{1}$, which is possible
since $\frac{b_{m}-s_{0}}{c_{m}-s_{0}}\rightarrow1.$ For sufficiently large
$m\geq m_{1},$

\qquad1. $\gamma_{m}$ and $\gamma$ belong to the same knot class and
$\gamma_{m}\in\mathcal{U}.$

\qquad2. $\left\{  s:\gamma_{m}(s)\neq\gamma(s)\right\}  \subset\lbrack
s_{0},b+A]$ which is disjoint from $I_{c}.$

\qquad3. $F_{k}(\gamma_{m})\geq F_{k}(\gamma),$ by construction of $f_{m}. $

\qquad4. $MDC(K_{m})\geq MDC(K)$, by Proposition 7(i).

\qquad5. $R_{O}(K_{m})\geq R_{O}(K),$ by Thickness Formula.

\qquad6. $\ell_{e}(\gamma_{m})\geq\ell_{e}(\gamma),$ since $\gamma$ is
relatively extremal and (1).

\qquad7. $\ell(\gamma_{m})\geq\ell(\gamma)$ by (5), (6) and the definition of
$\ell_{e}.$

\qquad8. $\ell(\gamma_{m})\leq\ell(\gamma)$ by construction of $f_{m}$ and
$\gamma_{m}.$

\qquad9. $f_{m}|[s_{0},c_{m}]$ and $\gamma|\left[  s_{0},b_{m}\right]  $ have
the same minimal length in

$\qquad\qquad\mathcal{C}(\gamma(s_{0}),\gamma(b_{m});\gamma^{\prime}%
(s_{0}),\gamma^{\prime}(b_{m});F_{k}(\gamma)^{-1}).$

\qquad10. $\gamma|[s_{0},b_{m}]$ is a CLC$(F_{k}(\gamma)^{-1})$-curve, by
Theorem 1 and $\kappa\gamma(s_{0})=0$.

We proved that\ if $\gamma|[s_{0},b+\varepsilon]$ is a $CLC(F_{k}(\gamma
)^{-1})$-curve (line segment followed by circular arc) such that
$0\leq\varepsilon<\pi F_{k}(\gamma)$ and $[s_{0},b+\varepsilon]\cap
I_{c}=\emptyset$, then there exists $\varepsilon_{1}=b_{m_{1}}-b>b_{0}%
-b=\varepsilon$ such that $\gamma|[s_{0},b+\varepsilon_{1}]$ is a
$CLC(F_{k}(\gamma)^{-1})$-curve. In fact, $\gamma|[s_{0},b+\varepsilon_{1}]$
must be one line segment followed by one circular arc by the definition of CLC
and the shape of $\gamma|[s_{0},b+\varepsilon]$.

Hence, $\varepsilon_{0}:=\max\left\{  \delta:\gamma|[s_{0},b+\delta]\text{ is
a }CLC(F_{k}(\gamma)^{-1})-\text{curve}\right\}  $ and $d=b+\varepsilon_{0}$
satisfies that $\varepsilon_{0}=\pi F_{k}(\gamma)$ or $d=b+\varepsilon_{0}\in
I_{c}.$ The proof is the same for the opposite direction before $a.$

i. Suppose that $\frac{1}{2}MDC(K)>R_{O}(K)=F_{k}(\gamma).$ One proceeds as in
proof of part (ii), omitting all conditions about avoiding $I_{c}$. Use
Proposition 8(ii), to obtain the line segment $\gamma|(a,b)$. Even though
$MDC(K_{m})\geq MDC(K)$ may not be valid by Proposition 7(i), $R_{O}%
(K_{m})\geq R_{O}(K)$ is valid by Proposition 7(ii). This shows that
$\gamma|[b,b+\pi F_{k}(\gamma)]$ is a $CLC(F_{k}(\gamma)^{-1})-$curve, even
passing through MDC-points. $\gamma(b)$ and $\gamma(b+\pi F_{k}(\gamma))$ is
an antipodal pair of a circle of radius $F_{k}(\gamma),$ forming a double
critical pair. This shows that $MDC(K)\leq\left\|  \gamma(b)-\gamma(b+\pi
F_{k}(\gamma)\right\|  =2F_{k}(\gamma)$ which is contrary to the hypothesis.
Consequently, the case of $\frac{1}{2}MDC(K)>R_{O}(K)=F_{k}(\gamma)$ with
$\exists s_{0}\in\mathbf{S}^{1},\kappa\gamma(s_{0})<\sup\kappa\gamma$ is vacuous.
\end{proof}

\begin{remark}
We do not know any example or the existence of any ideal knot $\gamma$ with
$R_{O}(K)=F_{k}(\gamma)<\frac{1}{2}MDC(K)$ and constant generalized curvature
$\kappa\gamma\equiv R_{O}(K)^{-1}.$
\end{remark}

\section{References}

[BS]\qquad G. Buck and J. Simon, \textit{Energy and lengths of knots,
}Lectures at Knots 96, 219-234.

[CE]\qquad J. Cheeger and D. G. Ebin, \textit{Comparison theorems in
Riemannian geometry, Vol 9, }North-Holland, Amsterdam, 1975.

[Di]\qquad Y. Diao, \textit{The lower bounds of the lengths of thick knots,
}preprint, August 2001.

[D]\qquad O. C. Durumeric, \textit{Thickness formula and }$C^{1}%
-$\textit{compactness of }$C^{1,1}$\textit{\ Riemannian submanifolds,
}preprint, January 2002.

[F]\qquad H. Federer, \textit{Geometric measure theory, }Springer, 1969.

[FHW]\qquad M. Freedman, Z.-X. He, and Z. Wang, \textit{Mobius energy of knots
and unknots}, Annals of Math, \textbf{139} (1994) 1-50.

[GM]\qquad O. Gonzales and H. Maddocks, \textit{Global curvature, thickness
and the ideal shapes of knots}, Proceedings of National Academy of Sciences,
\textbf{96 }(1999) 4769-4773.

[Ka]\qquad V. Katrich, J. Bendar, D. Michoud, R.G. Scharein, J. Dubochet and
A. Stasiak, \textit{Geometry and physics of knots}, Nature, \textbf{384}
(1996) 142-145.

[L]\qquad Litherland, \textit{Unbearable thickness of knots}, preprint.

[LSDR]\qquad A. Litherland, J Simon, O. Durumeric and E. Rawdon,
\textit{Thickness of knots}, Topology and its Applications, \textbf{91}(1999) 233-244.

[N]\qquad\ A. Nabutovsky, \textit{Non-recursive functions, knots ''with thick
ropes'' and self-clenching ''thick'' hyperspheres}, Communications on Pure and
Applied Mathematics\textit{, }\textbf{48} (1995) 381-428.

[RS]\qquad E. Rawdon and J. Simon, \textit{Mobius energy of thick knots},
Topology and its Applications, to appear.
\end{document}